\documentclass[a4paper,11pt]{article}
\usepackage{amsfonts,amssymb,latexsym,amsmath}
\usepackage{multicol}
\newcounter{num}[section]
\newcommand{\Num}{\refstepcounter{num}%
\textbf{\arabic{section}.\arabic{num}}}

\newcommand{\Theorem}{\textbf{Theorem~}}
\newcommand{\Proof}{\textbf{Proof}}
\newcommand{\Conj}{\textbf{Conjecture~}}
\newcommand{\Ex}{\textbf{Example}}
\newcommand{\Rem}{\textbf{Remark}}

\newcommand{\al}{\alpha}

\newcommand{\gx}{{\mathfrak g}}
\newcommand{\nx}{{\mathfrak n}}
\newcommand{\hx}{{\mathfrak h}}
\newcommand{\bx}{{\mathfrak b}}
\newcommand{\wdt}{{\dot{w}}}
\newcommand{\NH}{{\mathrm{Norm}(H)}}
\newcommand{\Ad}{{\mathrm{Ad}}}
\newcommand{\Oc}{{\cal O}}
\newcommand{\Ic}{{\cal I}}

\renewcommand{\leq}{\leqslant}

\begin{document}
\Large

\title{Tangent Cones of Schubert Varieties  for  $A_n$ of lower rank}
\author{A.N.Panov
\thanks{First author is supported by   RFBR-grant
09-01-00058}\and D.Yu.Eliseev}

\date{}

\maketitle

\section{Main conjectures}

Calculation of tangent cones for Schubert varieties at the origin
point is an interesting and extremely difficult problem. One of the
reasons is that   known methods of  calculation of tangent cones
 base on determination of Gr\"obner basis in the defining ideal of Schubert variety
 (more precisely, of its affine part). Even in the smooth case (when
 the tangent cone coincides with the tangent space) we have got an
 interesting and nontrivial theory (see \cite{L}).

In this paper we calculate tangent cones for series  $A_n$, where
$1\leq n\leq 4$, we formulate  the conjectures on structure of
tangent cones in general case. Determination of tangent cones is
important for classification of coadjoint orbits of maximal
unipotent subgroups (see \cite{K}), since each tangent cone is
 a subset stable with respect to the coadjoint representation.

Let  $G$  a semisimple  $K$-split algebraic group  over  a field $K$
of zero characteristic.   Lie algebra  $\gx$ of the group  $G$
admits decomposition  $\gx=\nx_-\oplus\hx\oplus \nx$, where $\hx$ is
the Cartan subalgebra, ~$\nx$ (resp. $\nx_-$) is a maximal nilpotent
subalgebra, spanned by the root vectors with positive (resp.
 оnegative) roots.  Denote by  $\bx=\hx\oplus \nx$ and $\bx_-=\hx\oplus
 \nx_-$. As usual, $H$,~$N$,~ $N_-$, ~$B$,~ $B_-$  are the
 corresponding subgroups of $G$. We denote by $\wdt$ an arbitrary chosen representative
 of the element $w$ of the Weyl group $W= \NH/H$.
  Using the Killing form we identify
   $\nx_-$ with the conjugate space $\nx^*$ of $\nx$.

 The group $G$ decomposes into the Bruhat classes
 $G= \bigcup_{w\in W} B\wdt
 B$.
 It implies that the flag variety  $X=G/B$ decomposes into the Schubert cells
  $X=\bigcup_{w\in W} X^0_w$, where $X^0_w=B\wdt
 B\bmod B$,~ $w\in W$. A closer  $X_w$ of the Schubert cell $X^0_w$ is called
 a Schubert variety.  Any Schubert variety contains
 an origin point  $p= B\bmod B$.

Denote by  $\Oc$ an affine open subset  $N_-B\bmod B$ in the flag
variety $X$. The set  $\Oc$  admits natural parametrization
$\exp(x)B\bmod B$, where $x\in \nx_-$ (for $\gx=A_n$ futher
$(1+x)B\bmod B$, ~ $x\in \nx_-$).

 The subset  $\Oc_w=\Oc\cap X_w$ is open in $X_w$ and closed in  $\Oc$.
 The origin point  $p$ belongs to  $\Oc_w$ and has zero
 coordinates following the chosen parametrization.
  We denote by  $C_w$ the tangent cone of  $X_w$ at the point  $p$
  (more precisely, it is a tangent cone of  $\Oc_w$ at zero point).

 By definition, for given closed subset  $M\subset K^n$, that contains
 the point $(0,\ldots,0)$, the tangent cone at zero is
 an annihilator of the ideal of lowest terms  $f_0$,
 where $f$ runs through the defining ideal   $I =I(M)$ (see
 \cite[chap.II,\S 1]{Shaf}
 or \cite[chap.9, \S 7]{KL}).

 The tangent cone  $C_w$ is contained in the tangent space
  $T_p(X)$ of the flag variety $X$ at the point  $p$. Identify
   $T_p(X)=\gx/\bx$
 with $\nx_-=\nx^*$. Since the subgroup  $B$ is a stabilizer
 of  $p$ in the group  $G$, the subgroup $B$ naturally acts in
  $T_p(X)=\nx^*$. This action coincides with the coadjoint action of  $B$
  on $\nx^*$. Any tangent cone $C_w$ is a  $\Ad^*$-invariant, closed
  with respect to  Zariski topology subset in  $\nx^*$.
  Notice first, that the the tangent cones may coincide for different elements
  $w\in W$. For instance, the tangent cones of Coxeter elements  coincide
  with $[\nx,\nx]^\perp$ (see \cite{K}).

  Well known that  $\dim X_w = l(w)$. Since the
dimension of an algebraic variety coincide with the dimension of its
tangent cone (see  theorem 8 in \cite[chap.9, \S 7]{KL}), we have
$\dim C_w = l(w)$.

The authors calculate the tangent cones $C_w$ for simple Lie
algebras of series $A_n$. The second author made the computer
program that calculate tangent cones. The algorithm of calculation
is presented in the second section.  This paper contains the results
of calculation for $ n\leq 4$, that are made by hand  and using
computer.
 Calculation of tangent cones for different examples
 provides  some general conjectures for an arbitrary that semisimple
 Lie algebras.

\Conj\Num. If  $C_{w_1} = C_{w_2}$, then the elements  $w_1$ and
$w_2$ are conjugate in the Weyl group. The converse statement is
false.

\Conj\Num. $C_w=C_{w^{-1}}$.

Let  $\gx_0$ be a semisimple regular subalgebra in  $\gx$ (a regular
subalgebra is a subalgebra stable with respect to adjoint action of
Cartan subgroup). The subalgebra  $\gx_0$ admits decomposition
$\gx=\nx_{0,-}\oplus\hx_0\oplus \nx_0$, where $\hx_0$, ~
$\nx_{0,-}$, ~ $\nx_0$  are the subalgebras in $\hx$,~ $\nx_-$,
~$\nx$ respectively. Identify the conjugate space  $\nx_0^*$ with
the subspace of $\nx^*$, that consists of all linear forms that
annihilate all roots vectors in  $\nx\setminus \nx_0$. The Weyl
group  $W_0$ of $\gx_0$ is a subgroup  of  $W$. For any  $w\in W_0$,
one can define both the tangent cone  $C_{w,0}$ in $\nx^*_0$ and
 the tangent cone  $C_w$ in $\nx^*$; note that $ C_{w,0}\subset
C_w$.

\Conj\Num. For any  $w\in W_0$ the closer  of $\Ad^*C_{w,0}$ is an
irreducible component in  $C_w$.

\Theorem\Num. The conjectures above are true for all Lie algebras
$A_n$ for  $n\leq 4$. \\
\Proof~~ follows from Tables 1-4.

\section{Algorithm of calculation of tangent cones and outcome}

For any root  $\gamma$ we denote: ~ $e_\gamma$ is a root vector, ~
 $x_\gamma(t)=\exp(te_\gamma)$, ~ $N_\gamma = \{x_\gamma(t): t\in K\}$, ~
 $N'_\gamma = \{x_\gamma(t): t\in
 K^*\}$.

 \Rem. If  $w=r_\al w'$ is a reduced decomposition, where  $\al$ is a simple root, then
$$ BwB = N_\al r_\al Bw'B \supset  N'_{\al}r_\al Bw'B =  N'_{-\al}Bw'B$$
(see \cite[\S 3, formulas  R1-R8 and Lemma 25]{Sh}). The subset
$N'_{-\al}Bw'B$ is dense in  $BwB$.

For any  $w\in W$ consider a reduced decomposition
$$w=r_{\al_1}\ldots r_{\al_l},$$ where $l=l(w)$ and $\al_1, \ldots,\al_l$
are  simple roots.

According to above Remark, the subset

$$ N'_{-\al_1} \cdots N'_{-\al_l}B$$  is dense in the Bruhat class  $BwB$.

Therefore the subset    $\Oc_w$ is a closer of the image of mapping
 $F: K^l\to \Oc$, where
$$F(t_1,\ldots,t_l) = x_{-\al_1}(t_1)\cdots x_{-\al_l}(t_l)\bmod
B.$$ The Elimination theory  (see ~\cite[\S 3]{KL}) provides the
method of construction of Gr\"obner basis of defining ideal $\Ic_w$
of the subset  $\Oc_w$. Further, using the standard procedure (see
Proposition 4(i) of chapter 9, \S 7  and Theorem 4 of chapter 8, \S
4 in the book  \cite{KL}) one can construct the Gr\"obner basis of
 the ideal $\Ic_{w,0}$, which annihilator coincides with $C_w$.
 Here is an example of calculation of
 tangent cone.\\
\Ex. $\gx = A_3$, ~$w=(13)(24)$. Identify  $\Oc$ with
$$N_-=\left\{\left(\begin{array}{cccc} 1&0&0&0\\
x_{21}&1&0&0\\
x_{31}&x_{21}&1&0\\
x_{41}&x_{42}&x_{43}&1\end{array}\right)\right\}.$$ To the reduced
decomposition  $w=(23)(12)(34)(23)$ we correspond the mapping
$F:K^4\to \Oc$, defined by formulas
$$x_{21}= t_2,\quad x_{31}=t_1t_2, \quad x_{41} = 0$$
$$x_{32}= t_1+t_4,\quad x_{42} = t_3t_4, \quad x_{43} = t_3.$$
Eliminating  $t_1$, $t_2$, $t_3$, $t_4$, we find the generators  $
x_{41}, x_{43}x_{31}+x_{42}x_{21}-x_{43}x_{32}x_{21}$ of $I(\Oc_w)$.
We obtain that the tangent cone is determined by the system of
equations  $x_{41}=0$,~ $x_{43}x_{31}+x_{42}x_{21}=0$.

We presents  the outcome of calculations of tangent cones for
$\gx=A_n$, ~$1\leq n\leq 4$.\\
{\bf Tangent cones for    $A_2$}.  The Weyl group coincides with
$S_3$. The equations that define the tangent cone in $$\nx^* = \left(\begin{array}{ccc}0&0&0\\
x_{21}&0&0\\x_{31}&x_{32}&0\end{array}\right),$$ are presented in
the following table.

\begin{center}
Table 1.
\end{center}

\begin{center}
\begin{tabular}{|c|c|}

\hline $w$&$C(w)$\\
\hline (13)& $\nx^*$\\
\hline (123),(132)& $x_{31}=0$\\
\hline (12)& $x_{31}=x_{32}=0$,~ \\
\hline (23)& $x_{31}=x_{21}=0$,~ \\
\hline e& $x_{31}=x_{32}=x_{21}=0$\\
\hline
\end{tabular}
\end{center}

{\bf Tangent cones for    $A_3$}.  The Weyl group coincides with
$S_4$. Introduce the notations $$D = \left|\begin{array}{cc}
x_{31}&x_{32}\\x_{41}&x_{42}\end{array}\right|, \quad P =
x_{43}x_{31}+x_{42}x_{21}.$$ The equations that define the tangent cone in
 $$\nx^* = \left(\begin{array}{cccc}0&0&0&0\\
x_{21}&0&0&0\\x_{31}&x_{32}&0&0\\x_{41}&x_{42}&x_{43}&0\\
\end{array}\right)$$ are presented in
the following table.

\begin{center}
Table 2.
\end{center}

\begin{center}
\begin{tabular}{|c|c|}

\hline $w$&$C(w)$\\
\hline (14)(23)& $\nx*$\\
\hline (14)& $ D=0$\\
\hline (1324), (1423)& $x_{41}=0$,~ \\
\hline (13)(24)& $x_{41}=0,~  P= 0$\\
\hline (134), (143)& $x_{41}=x_{42}=0$ \\
\hline (13)& $x_{41}= x_{42}=x_{43}=0$\\
\hline (124), (142)& $x_{41}= x_{31}= 0$\\
\hline (24)& $x_{41}= x_{31}=x_{21}=0$\\
\hline (1234), (1243), (1342), (1432)& $x_{41}= x_{31}= x_{42}=0$\\
\hline (234),  (243)& $x_{41}= x_{31}= x_{21}=x_{42}=0$\\
\hline (12)(34)& $x_{41}= x_{31}= x_{42}=x_{32}=0$\\
\hline (123),  (132)& $x_{41}= x_{31}=x_{42}=x_{43}=0$\\
\hline (12)& $x_{41}= x_{31}=x_{42}=x_{32}= x_{43}=0$\\
\hline (23)& $x_{41}= x_{31}=x_{21}=x_{42}=x_{43}=0$\\
\hline (34)& $x_{41}= x_{31}=x_{21}=x_{42} = x_{32}=0$\\
\hline e& $x_{41}= x_{31}=x_{21}=x_{42}=x_{32}= x_{43}=0$\\
\hline

\end{tabular}
\end{center}
{\bf Tangent cones for  $A_4$}.  The Weyl group coincides with
$S_5$. The equations that define the tangent cone in
$$\nx^* = \left(\begin{array}{ccccc}
0&0&0&0&0\\
x_{21}&0&0&0&0\\x_{31}&x_{32}&0&0&0\\x_{41}&x_{42}&x_{43}&0&0\\x_{51}&x_{52}&x_{53}&x_{54}&0\\
\end{array}\right)$$ are presented in
Tables 3 and 4.
\newpage

\begin{center}
Table  3.
\end{center}

\begin{center}
\begin{tabular}{|c|c|}
\hline
$w$  & $C(w)$\\
\hline
$(13425),(15243)$&$x_{51}=x_{41}=0$\\
\hline
$(14235),(15324)$&$x_{52}=x_{51}=0$\\
\hline
$(14325),(15234)$&$x_{51}=0,x_{41}x_{52}=0$\\
\hline
$(12345),(15432),(12453),$&\\
$(12354),(13542),(15432),$ & $x_{31}=x_{41}=x_{42}=x_{51}=x_{52}=x_{53}=0$\\
$(12543),(13452),(14532),(13542)$&\\
\hline
$(12435),(15342),(12534),(14352)$&$x_{31}=x_{41}=x_{51}=0=x_{52}=0$\\
\hline
$(13245),(15423),(14523),(13254)$&$x_{41}=x_{51}=x_{52}=0=x_{53}=0$\\
\hline
$(13524),(14253)$&$x_{41}=x_{51}=x_{52}=0, x_{31}x_{53}=0$\\
\hline
$(1425),(1524)$ & $x_{54}=0$\\
\hline
$(1325),(1523)$&$x_{51}=x_{41}=0,x_{53}x_{42}-x_{52}x_{43}=0$\\
\hline
$(1534),(1435)$&$x_{51}=x_{52}=0,x_{42}x_{31}-x_{41}x_{32}=0$\\
\hline
$(1324),(1423)$&$x_{41}=x_{51}=x_{52}=x_{53}=x_{54}=0$ \\
\hline
$(2435),(2534)$ & $x_{21}=x_{31}=x_{41}=x_{51}=x_{52}=0$\\
\hline
$(1352),(1235),(1253),(1532)$ & $x_{31}=x_{41}=x_{42}=x_{51}=x_{52}=0$\\
\hline
$(1245),(1254),(1542),(1452)$ & $x_{31}=x_{41}=x_{51}=x_{52}=x_{53}=0$\\
\hline
$(1354),(1345),(1453),(1543)$ & $x_{41}=x_{42}=x_{51}=x_{52}=x_{53}=0$\\
\hline
$(2345),(2354),(2543),(2453)$ & $x_{21}=x_{31}=x_{41}=x_{42}=x_{51}=x_{52}=x_{53}=0$\\
\hline
$(1234),(1243),(1432),(1342)$ & $x_{31}=x_{41}=x_{42}=x_{51}=x_{52}=x_{53}=x_{54}=0$\\
\hline
$(135),(153)$&$x_{41}=x_{42}=x_{51}=x_{52}=0$\\
\hline
$(125),(152)$&$x_{31}=x_{41}=x_{51}=0, x_{53}x_{42}-x_{52}x_{43}=0$\\
\hline
$(145),(154)$&$x_{51}=x_{52}=x_{53}=0, x_{42}x_{31}-x_{41}x_{32}=0$\\
\hline
$(124),(142)$&$x_{31}=x_{41}=x_{51}=x_{52}=x_{53}=x_{54}=0$\\
\hline
$(134),(143)$&$x_{41}=x_{42}=x_{51}=x_{52}=x_{53}=x_{54}=0$\\
\hline
$(235),(253)$&$x_{21}=x_{31}=x_{41}=x_{42}=x_{51}=x_{52}=0$\\
\hline
$(245),(254)$&$x_{31}=x_{41}=x_{42}=x_{51}=x_{52}=x_{53}=0$\\
\hline
$(234),(243)$&$x_{21}=x_{31}=x_{41}=x_{42}=x_{51}=x_{52}=x_{53}=x_{54}=0$\\
\hline
$(123),(132)$&$x_{31}=x_{32}=x_{41}=x_{42}=x_{51}=x_{52}=x_{53}=x_{54}=0$\\
\hline
$(345),(354)$&$x_{21}=x_{31}=x_{32}=x_{41}=x_{42}=x_{51}=x_{52}=x_{53}=0$\\
\hline
$w$ & $C(w)$\\
\hline
$(15)(24)$ & $\mathfrak{n^*}$\\
\hline
$(15)(234),(15)(243)$&$x_{52}x_{41}-x_{51}x_{42}=0$\\
\hline
\end{tabular}
\end{center}

\newpage
\begin{center}
Table 4.
\end{center}

\begin{center}

\begin{tabular}{|c|c|}
\hline
$(14)(25)$ & $x_{51}=0, x_{54}x_{41}+x_{53}x_{31}+x_{52}x_{21}=0$\\
\hline
$(15)(23)$ & $x_{52}x_{41}-x_{51}x_{42}=0, x_{41}x_{53}-x_{51}x_{43}=0, x_{53}x_{42}-x_{52}x_{43}=0$\\
\hline
$(15)(34)$&$x_{52}x_{41}-x_{51}x_{42}=0, x_{52}x_{31}-x_{51}x_{32}=0, x_{42}x_{31}-x_{41}x_{32}=0$\\
\hline
$(15)$&$x_{42}x_{31}-x_{41}x_{32}=0,x_{52}x_{41}-x_{51}x_{42}=0,$\\
$ $& $x_{53}x_{42}-x_{52}x_{43}=0,x_{52}x_{31}-x_{51}x_{32}=0,x_{53}x_{41}-x_{51}x_{43}=0$ \\
\hline
$(34)(125),(34)(152)$&$x_{51}=x_{41}=x_{31}=0$\\
\hline
$(23)(154),(23)(145)$&$x_{51}=x_{52}=x_{53}=0$\\
\hline
$(25)(134),(25)(143)$&$x_{51}=x_{41}=0,x_{54}x_{41}+x_{53}x_{31}+x_{52}x_{21}=0$\\
\hline
$(14)(253),(14)(235)$&$x_{51}=x_{52}=0,x_{54}x_{41}+x_{53}x_{31}+x_{52}x_{21}=0$\\
\hline
$(24)(135), (24)(153)$&$x_{41}=x_{51}=x_{52}=0$\\
\hline
$(13)(25)$&$x_{21}x_{42}+x_{43}x_{31}=0, x_{31}x_{53}+x_{21}x_{52}=0, x_{53}x_{42}-x_{52}x_{43}=0$\\
\hline
$(14)(35)$&$x_{32}x_{53}+x_{42}x_{54}=0, x_{54}x_{41}+x_{53}x_{31}=0, x_{42}x_{31}-x_{41}x_{32}=0$\\
\hline
$(14)(23)$&$x_{51}=x_{52}=x_{53}=x_{54}=0$\\
\hline
$(25)(34)$&$x_{21}=x_{31}=x_{41}=x_{51}=0$\\
\hline
$(13)(245),(13)(254)$&$x_{41}=x_{51}=x_{52}=x_{53}=0, x_{31}x_{43}+x_{21}x_{42}=0$\\
\hline
$(35)(124),(35)(142)$&$x_{31}=x_{41}=x_{51}=x_{52}=0, x_{42}x_{54}+x_{32}x_{53}=0$\\
\hline
$(14)$&$x_{51}=x_{52}=x_{53}=x_{54}=0, x_{42}x_{31}-x_{41}x_{32}=0$\\
\hline
$(25)$&$x_{21}=x_{31}=x_{41}=x_{51}=0, x_{53}x_{42}-x_{52}x_{43}=0$\\
\hline

$(12)(35)$&$x_{31}=x_{32}=x_{41}=x_{42}=x_{51}=x_{52}=0$\\
\hline
$(13)(45)$&$x_{41}=x_{42}=x_{51}=x_{52}=x_{53}=x_{54}=0$\\
\hline
$(13)(24)$&$x_{41}=x_{51}=x_{52}=x_{53}=x_{54}=0,x_{31}x_{43}+x_{21}x_{42}=0$\\
\hline
$(24)(35)$&$x_{21}=x_{31}=x_{41}=x_{51}=x_{52}=0,x_{42}x_{54}+x_{32}x_{53}=0$\\
\hline
$(12)(345),(12)(354)$&$x_{31}=x_{32}=x_{41}=x_{42}=x_{51}=x_{52}=x_{53}=0$\\
\hline
$(45)(132),(45)(123)$&$x_{31}=x_{41}=x_{42}=x_{43}=x_{51}=x_{52}=x_{53}=0$\\
\hline
$(13)$&$x_{41}=x_{42}=x_{43}=x_{51}=x_{52}=x_{53}=x_{54}=0$\\
\hline
$(35)$&$x_{21}=x_{31}=x_{32}=x_{41}=x_{42}=x_{51}=x_{52}=0$\\
\hline
$(24)$&$x_{21}=x_{31}=x_{41}=x_{51}=x_{52}=x_{53}=x_{54}=0$\\
\hline
$(12)(45)$&$x_{31}=x_{32}=x_{41}=x_{42}=x_{43}=x_{51}=x_{52}=x_{53}=0$\\
\hline
$(12)(34)$&$x_{31}=x_{32}=x_{41}=x_{42}=x_{51}=x_{52}=x_{53}=x_{54}=0$\\
\hline
$(23)(45)$&$x_{21}=x_{31}=x_{41}=x_{42}=x_{43}=x_{51}=x_{52}=x_{53}=0$\\
\hline
$(12)$&$x_{31}=x_{32}=x_{41}=x_{42}=x_{43}=x_{51}=x_{52}=x_{53}=x_{54}=0$\\
\hline
$(23)$&$x_{21}=x_{31}=x_{41}=x_{42}=x_{43}=x_{51}=x_{52}=x_{53}=x_{54}=0$\\
\hline
$(34)$&$x_{21}=x_{31}=x_{32}=x_{41}=x_{42}=x_{51}=x_{52}=x_{53}=x_{54}=0$\\
\hline
$(e)$& 0\\
\hline
\end{tabular}

\end{center}

\end{document}